\documentclass[a4paper, 10pt]{article}
\usepackage[all]{xy}

\usepackage[dvipdfmx]{graphicx}

\usepackage{latexsym}
\usepackage{amsfonts} 
\usepackage[mathscr]{euscript}
\usepackage{amsmath,mathrsfs,bm,amssymb,color,amscd}

\title{\textsf{Stability of  good quantum numbers in ground states }}
\date{\empty}

\author{
Tadahiro Miyao \\
{\it Department of Mathematics, Hokkaido University}\\
{\it Sapporo, Japan}\\
e-mail:
 miyao@math.sci.hokudai.ac.jp,
}

\newcommand{\one}{1}

\newcommand{\h}{\mathfrak{H}}

\newcommand{\ex}{e}

\newcommand{\R}{\mathrm{ran}}

\newcommand{\la}{\langle}
\newcommand{\ra}{\rangle}
\newcommand{\Tr}{\mathrm{tr}}

\newcommand{\slim}{\mbox{$\mathrm{s}$-$\displaystyle\lim_{n\to\infty}$}}

\newcommand{\BbbR}{\mathbb{R}}
\newcommand{\BbbN}{\mathbb{N}}

\newcommand{\BbbC}{\mathbb{C}}

\newcommand{\vphi}{\varphi}

\newcommand{\Cone}{\mathfrak{P}}

\newcommand{\no}{\nonumber \\}

\begin{document}

\newtheorem{define}{Definition}[section]
\newtheorem{Thm}[define]{Theorem}
\newtheorem{Prop}[define]{Proposition}
\newtheorem{lemm}[define]{Lemma}
\newtheorem{rem}[define]{Remark}
\newtheorem{assum}{Condition}
\newtheorem{example}{Example}
\newtheorem{coro}[define]{Corollary}

\maketitle

\begin{abstract}
Let $H$ be a self-adjoint operator, bounded from below and let $O$ be a bounded self-adjoint operator
 with purely discrete spectrum. Suppose that  (i) $E(H)=\inf \mathrm{spec}(H)$ is a simple eigenvalue, and (ii)
 $H$ strongly commutes with $O$. Let $\psi_H$ be the eigenvector associated with $E(H)$. By the assumptions (i) and (ii), $\psi_H$ is an eigenvector of $O$: $O\psi_H=\mu(H)\psi_H$. 
In the context of quantum mechanics, $\mu(H)$ is called a good quantum number.
In this note,
 we examine the stability of $\mu(H)$ under perturbations of $H$ from a viewpoint of the order theory.
 In addition, we provide some applications of the theory to the study of ferromagnetism.
 \begin{flushleft}
{\bf Mathematics Subject Classification (2010).} 
\end{flushleft}
Primary:  47A63, 47B60,  47A75;
Secondary: 47A55
\begin{flushleft}
{\bf
Keywords. 
} 
\end{flushleft}
 Ground states; Eigenvalues;
Positivity improving resolvents; Posets;  Good quantum numbers; ferromagnetism.
\end{abstract}

\section{Introduction}
Let $\h$ be a complex Hilbert space and let $H$ be a self-adjoint operator on $\h$, bounded from below.
Suppose that $E(H)=\inf \mathrm{spec}(H)$ is a simple  eigenvalue, where $\mathrm{spec}(H)$ is spectrum of $H$.  The eigenvector associated with $E(H)$ is denoted by  $\psi_H$.
Let $O$ be a bounded self-adjoint operator with  purely discrete spectrum.
Assume  that $H$  strongly  commutes with $O$, that is, their spectral measures commute with each other.
Under this setting,  we readily see that $\psi_H$ is an eigenvector of $O$:
\begin{align}
O\psi_H=\mu(H)\psi_H.
\end{align}
In quantum mechanics, suppose that  a particular Hamiltonian $H$ and an operator $O$ with corresponding eigenvalues and eigenvectors are given. Then the eigenvalues are said to be  ``good quantum numbers" if every eigenvector remains an eigenvector of the same eigenvalue as time evolves, or $H$ strongly commutes with $O$. 
In this sense, the  eigenvalue $\mu(H)$  can be regarded as    a good quantum number.
\medskip

In this note, we will examine the stability of $\mu(H)$. To be precise, let $V$
 be a self-adjoint operator. 
 We will consider a perturbation of $H$ by $V$.
 For simplicity, we suppose  that $V$ is bounded.
 We continue to assume that  
 $E(H+V)$ is a simple eigenvalue of $H+V$.
Our main purpose  is to answer the following question.
\begin{center}
{\it
When does the equation  $\mu(H+V)=\mu(H)$ hold ?
}
\end{center}
In the rest of the present note, we will provide a framework which enables us to solve the above problem. 
Our novel idea for constructing the framework is to apply the positivity improvingness of the resolvent  of $H$. 
\medskip

Before we proceed, we  briefly explain the motivation behind the aforementioned problem.
An essence of the problem  originates from the study of ferromagnetism in many-electron systems;
Mathematical studies of ferromagnetism were initiated by Lieb \cite{Lieb}, Nagaoka \cite{Nagaoka} and Thouless \cite{Thouless}.
The origin of ferromagnetism still has been   mystery and are actively examined even today, see e.g.,  \cite{KT, KSV, Tasaki3}.
In \cite{Miyao5,Miyao6,Miyao7,Miyao8}, Miyao examined the stability of ferromagnetism in many-electron systems. In particular, he gave 
a model independent framework which describes various stability results concerning
 ferromagnetism in  the Hubbard model \cite{Miyao8}.
Remark that in concrete applications to many-electron system, $H$ corresponds to the Hamiltonian and $O$ corresponds to the total spin operator. 
In the present note, we focus our attention on a  mathematical aspect of the theory established in \cite{Miyao8}. We will find  that its structure 
is well  decribed by the order theory.
\medskip
 
The rest of the present note is organized as follows.
In Section \ref{Sec2}, we introduce some basic notions to state our main result.
In particular, we focus our attention on the study of  positivity improving resolvents.
In Section \ref{Sec3}, we state our main result; we provide a novel framework which solve the stability problem stated in this introduction.  Section \ref{Sec4} is devoted to give an example. This example suggests that our  framework  contains rich mathematical strucutures. In Appendices \ref{AppA}
 and \ref{AppB}, we prove some operator inequalities which are useful in the main sections.

\begin{flushleft}
{\bf Acknowledgements}
\end{flushleft}
The author  was partially supported by KAKENHI 18K03315.

\section{Preliminary}\label{Sec2}
\subsection{Basic definitions}
Let $\h$ be a complex Hilbert space.
By a {\it convex  cone}, we understand a closed convex set  $\Cone\subset \h$
such that $t\Cone \subseteq \Cone$ for all $t\ge 0$ and $\Cone\cap (-\Cone)=\{0\}$.
The {\it dual cone   of} $\Cone$ is defined by 
$
\Cone^{\dagger}=\{\eta\in \h\, |\, \la \eta|\xi\ra\ge 0\ \forall \xi\in
\Cone \}.
$
 We say that $\Cone$ is {\it self-dual} if 
$
\Cone=\Cone^{\dagger}. $
 In what follows, we always assume that $\Cone$ is self-dual.

 A vector $\xi$ is said to be  {\it positive w.r.t. $\Cone$} if $\xi\in
 \Cone$.  We write this as $\xi\ge 0$  w.r.t. $\Cone$.
  A vector $\eta\in \Cone$ is called {\it strictly positive
w.r.t. $\Cone$} whenever $\la \xi| \eta\ra>0$ for all $\xi\in
\Cone\backslash \{0\}$. We write this as $\eta>0 $
w.r.t. $\Cone$.

We denote by  $\mathscr{B}(\h)$  the set of all bounded linear operators on
$\h$.
\begin{define}{\rm 

Let $A\in \mathscr{B}(\h)$. 
\begin{itemize}
 \item If $A \Cone\subseteq \Cone,$\footnote{
For each subset $\mathfrak{C}\subseteq \h$, $A\mathfrak{C}$ is
	     defined by $A\mathfrak{C}=\{Ax\, |\, x\in \mathfrak{C}\}$.
} we then 
write  this as  $A \unrhd 0$ w.r.t. $\Cone$. In
	     this case, we say that {\it $A$ preserves the
positivity w.r.t. $\Cone$.}  
\item The set of all positivity preserving operators w.r.t. $\Cone$ is denoted by $\mathscr{P}(\Cone)$.
$\mathscr{P}(\Cone)$ is a weakly closed convex cone in $\mathscr{B}(\h)$.
\item 
We write  $A\rhd 0$ w.r.t. $\Cone$, if  $A\xi >0$ w.r.t. $\Cone$ for all $\xi\in
\Cone \backslash \{0\}$. 
 In this case, we say that {\it $A$ improves the
positivity w.r.t. $\Cone$.} 
\end{itemize}
} 
\end{define} 
Let $\mathscr{L}^2(\h)$ be the set of all Hilbert-Schmidt class operators on  $\h$. 
It is well-known that $\mathscr{L}^2(\h)$ becomes a Hilbert space by introducing the norm 
$\la A|B\ra_{\mathscr{L}^2}=\Tr[A^*B],\ \ A, B\in \mathscr{L}^2(\h)$.

Let $\mathfrak{S}(\h)$ be the set of  all density matries on $\h$.\footnote{
We say that a linear operator $\varrho$ on $\h$ is 
a  {\it density matrix}, if it is     positive  and  trace class operator such that $\Tr[\varrho]=1$.}
For given normalized vector $\psi$, we set $\varrho_{\psi}=|\psi\ra\la\psi| $. Trivially, we have $\varrho_{\psi} \in \mathfrak{S}(\h)$. If $\psi\ge 0$ w.r.t. $\Cone$, then $
\varrho_{\psi} \unrhd 0 $ w.r.t. $\Cone$, while, if $\psi>0$ w.r.t. $\Cone$, then $\varrho_{\psi} \rhd 0$
w.r.t. $\Cone$.

In the present note, we will examine self-adjoint operators  satisfying the following conditions:
\begin{itemize}
\item[1.] $H$ is self-adjoint and bounded from below;
\item[2.]  $H$ has purely discrete spectrum; 
\item[3.] $(H+s)^{-1} \unrhd 0$ w.r.t. $\Cone$ for all $s > -E(H)$, where $E(H)=\inf \mathrm{spec}(H)$.

\end{itemize}
We denote by $\mathscr{A}_{\Cone}$ the set of all operators satisfying the conditions 1.-3. above.

\begin{Prop}\label{ConvexConeOp}
Let $H, H'\in \mathscr{A}_{\Cone}$. If $sH+tH'$ is essentially self-adjoint for some  $s>0$ and $t>0$, then 
$\overline{sH+tH'} \in \mathscr{A}_{\Cone}$. In particular, 
$\mathscr{A}_{\Cone} $ is a convex cone.
\end{Prop}
{\it Proof.} 
By Proposition \ref{Equiv}, $e^{-\beta sH} \unrhd 0$ and $e^{-\beta t H'} \unrhd 0$ w.r.t. $\Cone$ for all $\beta \ge 0$.
By the Trotter product  formula \cite[Theorem S.20]{ReSi1}, we have
\begin{align}
e^{
-\beta(
\overline{sH+tH'})
}
=\slim \Big(
e^{-\beta sH/n} e^{-\beta tH'/n}
\Big)^n, \label{Trotter}
\end{align}
where $\slim$ indicates the strong limit.
Because $e^{-\beta sH/n} \unrhd 0$ and $e^{-\beta tH'/n} \unrhd 0$ w.r.t. $\Cone$, we see that 
$
\big(
e^{-\beta sH/n} e^{-\beta tH'/n}
\big)^n\unrhd 0
$ w.r.t. $\Cone$ for all $\beta \ge 0$ and $n\in \BbbN$. Thus, the right hand side of (\ref{Trotter})
 preserves the positivity w.r.t. $\Cone$ for all $\beta \ge 0$. 
 By applying Proposition \ref{Equiv} again, we obtain that $\overline{sH+tH'}\in \mathscr{A}_{\Cone}$.
 $\Box$
\medskip

Let $\mathscr{A}^+_{\Cone}$ is the set of all  self-adjoint operators  satisfying 1., 2. and 3'. below:
\begin{itemize}
\item[3'.] $(H+s)^{-1}\rhd 0$ w.r.t. $\Cone$ for all $s> -E(H)$.
\end{itemize}
\begin{rem}
{\rm 
If $H\in \mathscr{A}^+_{\Cone}$, then $E(H)$ is a simple eigenvalue with strictly positive eigenvector 
by Theorem \ref{Uni}.
}
\end{rem}

\subsection{Propagation of positivity}
Let $\h_1$ and $\h_2$ be complex Hilbert spaces,  and let $\Cone_1$ and $\Cone_2$ be self-dual cones in $\h_1$ and $\h_2$, respectively.
Suppose that $\h_1$ is a closed subspace of $\h_2$.
The orthogonal projection from $\h_2$ to $\h_1$ is denoted by $\pi_{1, 2}$.
We say that the {\it  positivity  is inherited from  $\Cone_1$ to $\Cone_2$,  } if
the following are satisfied:
\begin{itemize}
\item[1.] $\Cone_1=\pi_{1, 2}\Cone_2$;
\item[2.] $\pi_{1, 2} \unrhd 0$ w.r.t. $\Cone_2$.
\end{itemize}
In this case, we write $\Cone_1\dashrightarrow \Cone_2$. As we will see, this binary relation defines a partial order\footnote{Readers are referred to \cite{Roman}
for partial orders.}.
The {\it conditional expectation } $\mathscr{E}_{1, 2}: \mathscr{B}(\h_2)\to \mathscr{B}(\h_2)$ is defined by 
\begin{align}
\mathscr{E}_{1, 2}(A) =\pi_{1, 2} A\pi_{1, 2}+\pi_{1, 2}^{\perp} A\pi_{1, 2}^{\perp},\ \ A\in \mathscr{B}(\h_2).
\end{align}
The following proposition is often useful:
\begin{Prop}\label{PPtoPP}
$
\mathscr{E}_{1, 2}\big(\pi_{1, 2}\mathscr{P}(\Cone_2)\pi_{1, 2}\big) = \pi_{1, 2}\mathscr{E}_{1, 2}\big(\mathscr{P}(\Cone_2)\big) \pi_{1, 2}= \mathscr{P}(\Cone_1).
$
\end{Prop}
{\it Proof.} Let $A\in \mathscr{P}(\Cone_2)$.
We have 
\begin{align}
\pi_{1, 2} A\pi_{1, 2} \Cone_1=\pi_{1, 2} A\pi_{1, 2} \Cone_2\subseteq \pi_{1, 2} A\Cone_2\subseteq \pi_{1, 2} \Cone_2=\Cone_1.
\end{align}
Therefore, we have $\mathscr{E}_{1, 2}(\pi_{1 , 2}\mathscr{P}(\Cone_2)\pi_{1, 2}) \subseteq  \mathscr{P}(\Cone_1)$.
Conversely, suppose that $A\in \mathscr{P}(\Cone_1)$. Then, corresponding to the decomposition $
\h_2=\h_1\oplus \h_1^{\perp}$, we have
 $A\oplus 0 \in \mathscr{P}(\Cone_2)$.
Indeed, for all $x, y\in \Cone_2$, we have
$\la x|A\oplus 0 y\ra=\la \pi_{1,2} x|A\pi_{1, 2}y\ra \ge 0$. Because $
\mathscr{E}_{1, 2} (A\oplus 0)=A\oplus 0
$, we obtain that $\pi_{1,2 }\mathscr{E}_{1, 2}(\mathscr{P}(\Cone_2)) \pi_{1, 2} \supseteq  \mathscr{P}(\Cone_1)$. $\Box$

\begin{define}
{\rm 
Let $H_1\in \mathscr{A}_{\Cone_1}$ and $H_2\in \mathscr{A}_{\Cone_2}$ be self-adjoint operators bounded from below.
If $\Cone_1 \dashrightarrow \Cone_2$ is satisfied, then we say that 
{\it the  $\Cone_2$-positivity  of $H_2$ is inherited from the  $\Cone_1$-positivity  of $H_1$ } and write this as  
$(H_1, \Cone_1) \dashrightarrow (H_2, \Cone_2)$.
}
\end{define}

\begin{Prop}
Suppose that $(H_1, \Cone_1) \dashrightarrow (H_2, \Cone_2)$.
Suppose that $E(H_1)$ and $E(H_2)$ are simple eigenvalues.
Then the corresponding eigenvectors, say $\psi_{H_1}$ and $\psi_{H_2}$,  are positive, namely, 
$\psi_{H_1} \ge 0$ w.r.t.  $\Cone_1$, and $\psi_{H_2} \ge 0$  w.r.t. 
 $\Cone_2$. Equivalently, we have $\varrho_{\psi_{H_1}}\unrhd 0$ w.r.t. $\Cone_1$ and $\varrho_{\psi_{H_2}} \unrhd 0$ w.r.t. $\Cone_2$.
Moreover,    $ \pi_{1, 2}\mathscr{E}_{1, 2}(\varrho_{\psi_{H_2}}) \pi_{1, 2} \unrhd 0$ w.r.t. $\Cone_1$ and $\la \psi_{H_1}|\pi_{1, 2} \psi_{H_2} \ra \ge 0$  hold.
\end{Prop}
{\it Proof.} By Proposition \ref{GSP}, we have $\psi_{H_1}\ge 0$ w.r.t. $\Cone_1$ and $\psi_{H_2}\ge 0$ w.r.t. $\Cone_2$, respectively.
From the property  $\Cone_1=\pi_{1, 2}\Cone_2$, it holds that $\pi_{1, 2}\psi_{H_2}\ge 0$ w.r.t. $\Cone_1$.
Thus, we obtain that  $\la \psi_{H_1}|\pi_{1, 2} \psi_{H_2} \ra \ge 0$. 
Because $\varrho_{\psi_{H_2}} \unrhd 0$ w.r.t. $\Cone_2$, we obtain $ \pi_{1, 2}\mathscr{E}_{1, 2}(\varrho_{\psi_{H_2}}) \pi_{1, 2}\unrhd 0$ by Proposition \ref{PPtoPP}.
$\Box$
\medskip

Let $\{\h_n\}_{n=1}^{\infty}$ be a sequence of Hilbert spaces, and let $\{\Cone_n\}_{n=1}^{\infty}$ 
and  $\{\Cone_n'\}_{n=1}^{\infty }$ be 
 sequences of self-dual cones such that
\begin{itemize}
\item[(i)] $\h_n$ is a closed subspace of $\h_{n+1}$;
\item[(ii)] $\Cone_n$ and $\Cone_n'$  are   self-dual cones in $\h_n$.
\end{itemize}
Assume that $\mathscr{A}_{\Cone_n} \cap \mathscr{A}_{\Cone_n'} \neq \{0\}$ for all $n\in \BbbN$.
Suppose that a sequence of semibounded self-adjoint operators $\{H_n\}_{n=1}^{\infty}$ satisfies the 
following:
\begin{align}
 (H_1, \Cone_1) &\dashrightarrow (H_2, \Cone_2'),  (H_2, \Cone_2) \dashrightarrow (H_3, \Cone_3'), 
  (H_3, \Cone_3)  \dashrightarrow (H_4, \Cone_4') , (H_4, \Cone_4) \dashrightarrow \no
  \dashrightarrow \cdots &\dashrightarrow (H_n, \Cone_n'),  (H_n, \Cone_n) \dashrightarrow (H_{n+1}, \Cone_{n+1}'), (H_{n+1}, \Cone_{n+1}) \dashrightarrow \cdots.
\end{align}
By definition, we have $\la \psi_{H_1}|\pi_{1, 2} \psi_{H_2}\ra \la \psi_{H_2}|\pi_{2, 3} \psi_{H_3}\ra
\cdots \la \psi_{H_n}|\pi_{n, n+1} \psi_{H_{n+1}}\ra \ge 0$
and $\pi_{j, j+1} \mathscr{E}_{j, j+1}(\varrho_{\psi_{H_{j+1}}}) \pi_{j, j+1} \unrhd 0$
w.r.t. $\Cone_j$, where $\pi_{j, j+1}$ is the orthogonal projection from $\h_{j+1}$ onto $\h_j$.
In this sense, the positivity of $\psi_{H_1}$ is  propageted  to $\psi_{H_{n}}$ for every $n\in \BbbN$.

\subsection{Propagation of  strict positivity}
\begin{define}
{\rm 
Let $H_1\in \mathscr{A}_{\Cone_1}^+$ and $H_2\in \mathscr{A}_{\Cone_2}^+$.
If $\Cone_1 \dashrightarrow \Cone_2$ is satisfied, then we say that {\it  the   strict 
 $\Cone_2$-positivity of $H_2$ is inherited from the  stirct $\Cone_1$-positivity of $H_1$ } and write this as  
$(H_1, \Cone_1) \to  (H_2, \Cone_2)$. By  definition, we readily confirm that if $
(H_1, \Cone_1) \to (H_2, \Cone_2)
$, then we have $ (H_1, \Cone_1) \dashrightarrow (H_2, \Cone_2)$.
}
\end{define}

Let $H_1\in \mathscr{A}_{\Cone_1}^+$ and $H_2\in \mathscr{A}_{\Cone}^{+}$. As before,
 the ground state of $H_j$ is denoted by $\psi_{H_j},\ j=1, 2$. 

\begin{Thm}\label{StrictP}
 If $
(H_1, \Cone_1) \to (H_2, \Cone_2)$, then
  $\la \psi_{H_1}|\pi_{1, 2} \psi_{H_2}\ra>0$.
   Equivalently, 
   we have  $
   \big\la \varrho_{\psi_{H_1}}\big|\mathscr{E}_{1, 2}(\varrho_{\psi_{H_2}}) \big\ra_{\mathscr{L}^2}>0
   $. In addition, $\pi_{1, 2} \mathscr{E}_{1, 2}(\varrho_{\psi_{H_2}}) \pi_{2, 1} \rhd 0$ w.r.t. $\Cone_1$.
\end{Thm}

To prove Theorem \ref{StrictP}, we begin with the following lemma:

\begin{lemm}\label{VecP2}
Let $A\in \mathscr{B}(\h)$ with $A\neq 0$. Assume that  $u>0$ w.r.t. $\Cone$.
If $A\unrhd 0$ w.r.t. $\Cone$, then $Au\neq 0$.
\end{lemm} 
{\it Proof.}
First, we prove the following claim:
 Let $A\in \mathscr{B}(\h)$. If $Au= 0$ for all $u\in \Cone$, then $A=0$.
 Indeed, by Proposition \ref{Haagerup}, each $u\in \h$ can be written as 
$
u=v_1-v_2+i(w_1-w_2)
$, where $v_1, v_2, w_1, w_2\in \Cone$ such that $\la v_1|v_2\ra=0$ and
$\la w_1|w_2\ra=0$. Thus, the assumption implies that $Au=0$
 for {\it all} $u\in \h$.

Assume that $Au=0$. Then,  $\la v|Au\ra=0$ for all $v\in \Cone$,
implying that  $\la A^*v|u\ra=0$.
Since $u>0$  and $A^* v\ge 0$ w.r.t. $\Cone$, we conclude that 
$A^*v$ must be zero. Because $v$ is arbitrary, 
$A^*=0$  by the above claim. 
This contradicts with the assumption $A\neq 0$. 
$\Box$

\begin{flushleft}
{\it Proof of Theorem \ref{StrictP}}
\end{flushleft}
Note that $\psi_{H_1} >0$ w.r.t. $\Cone_1$ and $\psi_{H_2}>0$ w.r.t. $\Cone_2$, respectively.
Because $\Cone_1=\pi_{1, 2}\Cone_2$ and  $\pi_{1, 2} \unrhd 0$ w.r.t. $\Cone_2$, we obtain that $\pi_{1, 2} \psi_{H_2} \ge 0$ w.r.t. $\Cone_1$ and $\pi_{1, 2}\psi_{H_2}\neq 0$ by Lemma \ref{VecP2}. Since $\psi_{H_1} >0$ w.r.t. $\Cone_1$, we conclude that 
$\la \psi_{H_1}|\pi_{1, 2} \psi_{H_2}\ra >0$.

For each $x, y\in \Cone \backslash \{0\}$, we have
$\la x |\mathscr{E}_{1, 2}(\varrho_{\psi_{H_2}}) y\ra=\la x\oplus 0 | \psi_{H_2}\ra \la \psi_{H_2}|y\oplus 0\ra>0$, which implies that $\pi_{1, 2} \mathscr{E}(\varrho_{\psi_{H_2}}) \pi_{1, 2} \rhd 0$ w.r.t. $\Cone_1$.
  $\Box$
\medskip

As before, assume that $\mathscr{A}_{\Cone_n} \cap \mathscr{A}_{\Cone_n'} \neq \{0\}$ for all $n=1, \dots, N$.
Suppose that a sequence of semibounded self-adjoint operators $\{H_n\}_{n=1}^{N}$ satisfies the 
following:
\begin{align}
 (H_1, \Cone_1) &\to  (H_2, \Cone_2'),  (H_2, \Cone_2) \to (H_3, \Cone_3'), 
  (H_3, \Cone_3)  \to (H_4, \Cone_4') , (H_4, \Cone_4) \to  \no
  \to  \cdots &\to  (H_{N-1}, \Cone_{N-1}'),  (H_{N-1}, \Cone_{N-1}) \to (H_{N}, \Cone_{N}').\label{ChainE}
\end{align}
Then,  we have $\la \psi_{H_1}|\pi_{1, 2} \psi_{H_2}\ra \la \psi_{H_2}|\pi_{2, 3} \psi_{H_3}\ra
\cdots \la \psi_{H_{N-1}}|\pi_{{N-1}, N} \psi_{H_{N}}\ra > 0$ for each $N\in \BbbN$ (equivalently, 
$
\big\la \varrho_{\psi_{H_1}}\big|\mathscr{E}_{1, 2}(\varrho_{\psi_{H_2}}) \big\ra_{\mathscr{L}^2}
\cdots \big\la \varrho_{\psi_{H_{N-1}}}\big|\mathscr{E}_{N-1, N}(\varrho_{\psi_{H_N}}) \big\ra_{\mathscr{L}^2}
>0$)
    by Theorem \ref{StrictP}. Furthermore, $
    \pi_{j, j+1} \mathscr{E}_{j, j+1}(\varrho_{H_{j+1}}) \pi_{j, j+1} \rhd 0
    $ w.r.t. $\Cone_{j}$ for all $j\in \BbbN$.
In this sense,  a strict positivity of $\psi_{H_1}$ is  propageted to $\psi_{H_N}$. 
As we will see in the following sections, this property is  important to examine the stability of the good quantum numbers. 

\begin{define}
{\rm 
We say that $(H_1, \Cone_1)$ and $(H_N, \Cone_N')$ are connected by  
the sequences $\{(H_j, \Cone_j)\}_{j=1}^{N-1}$ and $\{(H_j, \Cone_j')\}_{j=2}^N$ if 
 (\ref{ChainE}) holds.
We simply express   this  as $H_1  \to  H_N$.
}
\end{define}

For a given Hilbert space $\h_*$, let $\mathbb{H}_{\h_*}$ be the set of all Hilbert spaces  containing 
$\h_*$ as a closed subspace. Let $\mathscr{P}_{\h_*, 0}$ be the set of self-adjoint operators defined by 
\begin{align}
\mathscr{P}_{\h_*, 0}=\bigcup_{\h\in \mathbb{H}_{\h_*}} \bigcup_{\Cone\subset \h} \mathscr{A}_{\Cone}^+,
\end{align}
where the union $\bigcup_{\Cone\subset \h}$ runs  over   all self-dual cones in $\h$.

\begin{Prop}\label{PreO}
The binary relation `` $\to$'' is a preoder on $\mathscr{P}_{\h_*, 0}$. Namely, we have the following:
\begin{itemize}
\item[{\rm (i)}] $H\to H$; 
\item[{\rm (ii)}] $H \to H',\ \  H' \to  H'' \Longrightarrow H \to
 H'' $.
\end{itemize} 
\end{Prop}
{\it Proof.} (i) Because $H\in \mathscr{P}_{\h_*, 0}$, there is a self-dual cone $\Cone$ such that 
$H\in \mathscr{A}_{\Cone}^+$. Then we can  readily check that $(H, \Cone)\to (H, \Cone)$.

(ii) By definition, $H$ and $H'$ are connected by seqences $P=\{(H_j, \Cone_j)\}_{j=1}^{N-1}$ and $P'=\{H_j, \Cone_j'\}_{j=2}^{N}$
 with $H_1=H$ and $H_N=H'$. Also $H'$ and $H''$ are connected by sequences 
 $Q=\{(K_j, \mathfrak{Q}_j)\}_{j=1}^{M-1}$ and $Q'=\{(K_j, \mathfrak{Q}_j')\}_{j=2}^{M}$
 with $K_1=H'$ and $K_M=H''$.
Now we define new sequences $R$ and $R'$ by $
R=P\cup Q
$ and $R'=P\cup Q'$, then $H$ and $H''$ are connected by $R$ and $R'$. $\Box$

\begin{define}\label{HEquiv}
{\rm 
Let $H_1, H_2\in \mathscr{P}_{\h_*, 0}$.
If $H_1\to H_2$ and $H_2\to H_1$, then we write this as  $H_1\equiv H_2$.
The binary relation \lq\lq{} $\equiv$ \rq\rq{} is an equivalence relation on $\mathscr{P}_{\h_*, 0}$.
Let $\mathscr{P}_{\h_*, 0}$ be the set of equivalence classes:  $\mathscr{P}_{\h_*}=\mathscr{P}_{\h_*, 0} / \equiv$.
The equivalence class containing $H$ is denoted  by $[H]$.
The binary relation \lq\lq{}$\to $\rq\rq{} on $\mathscr{P}_{\h_*}$ is naturally  defined by 
$
[H_1] \to  [H_2]\ \ \ \mbox{if}\ \ H_1\to  H_2.
$
This is a {\it partial order} on $\mathscr{P}_{\h_*}$; namely, we have, by Proposition \ref{PreO},
\begin{itemize}
\item[(i)] $[H] \to [H]$;
\item[(ii)] $[H_1] \to [H_2],\ [H_2] \to [H_1]\Longrightarrow [H_1] =[H_2]$; 
\item[ (iii)] $[H_1]\to [H_2],\   [H_2] \to  [H_3] \Longrightarrow [H_1]\to
 [H_3]$. 
\end{itemize}
In what follows, we abbreviate $[H]\to [H']$ to $H\to H'$ if no confusion arises.
}
\end{define}

\section{Stability of good quantum numbers in ground states}\label{Sec3}

\subsection{Main result}

Let $O\in \mathscr{B}(\mathfrak{H}_*)$ be self-adjoint.
In what follows, we always assume that $O$ has  purely discrete spectrum.
In this section, we will explore the following class of self-adjoint operators:
\begin{align}
\mathscr{P}_{\h_*, 0}(O)=\{H\in \mathscr{P}_{\h_*, 0}\, |\, \mbox{$e^{isO} e^{itH}=e^{itH}e^{isO}$ for all $s, t\in \BbbR$}\}.
\end{align}
For each $\h\in \mathbb{H}_{\h_*}$, $O$ can be naturally extended to a self-adjoint operator on  $\h$.\footnote{
Indeed, corresponding to the decomposition $\h=\h_*\oplus \h_*^{\perp}$,
the natural extension of $O$ to $\h$ is defined by $O\oplus 0$. We denote by $O$ this natural extension 
if no confusion arises.
} The natural extension is also denoted  by the same symbol $O$.

Note that the preorder ``$\to$'' can be  defined on   $\mathscr{P}_{\h_*, 0}(O)$ as well. 
As before, we set $\mathscr{P}_{\h_*}(O)=\mathscr{P}_{\h_*, 0}(O) / \equiv$.
Then the preoder ``$\to$'' can be also lifted  up to a  partial order.
We identify  the equivalence class $[H]\in \mathscr{P}_{\h_*}(O)$ with $H$ if no confusion occurs.

\begin{Prop}\label{mumu}
Let $H, K\in \mathscr{P}_{\h_*}(O)$.
If $H\to K$, then $\mu(H)=\mu(K)$.
\end{Prop}
{\it Proof.}  Suppose that $H\in \mathscr{A}_{\Cone}^+$ and $K\in \mathscr{A}_{\mathfrak{Q}}^+$.
There exist   sequences $\{(H_j, \Cone_j)\}_{j=1}^{N-1}$ and $\{(H_j, \Cone_j')\}_{j=2}^N$ satisfying (\ref{ChainE}) with $H_1=H,\ \Cone_1=\Cone,\ H_N=K$ and $\Cone_N'=\mathfrak{Q}$.
Let $\pi_{j, j+1}$ be the orthogonal projection from $\h_{j+1}$ onto $\h_j$.
Let $\mathscr{E}_{j, j+1}$ be the corresponding  conditional expectation. Because $O\in \mathscr{B}(\h_*)$,
we see that $\mathscr{E}_{j, j+1}(O)=O$, which implies that 
\begin{align}
\mu(H_j) \la \psi_{H_j}|\pi_{j, j+1}\psi_{H_{j+1}}\ra&=\la O \psi_{H_j}|\pi_{j, j+1}\psi_{H_{j+1}}\ra\no
&=\la \psi_{H_j}|\mathscr{E}_{j, j+1}(O)\psi_{H_{j+1}}\ra\no
&=\la \psi_{H_j}|\pi_{j, j+1}O\psi_{H_{j+1}}\ra\no
&=\mu(H_{j+1}) \la \psi_{H_j}|\pi_{j, j+1}\psi_{H_{j+1}}\ra.
\end{align}
By applying Theorem \ref{StrictP}, we obtain that
$\mu(H_j)=\mu(H_{j+1})$. Repeating this argument several times, we arrive at  
$\mu(H)=\mu(H_1)=\mu(H_2)=\cdots=\mu(H_N)=\mu(K)$. $\Box$

\begin{define}\label{DefUniO2}
{\rm
Let $H_*\in \mathscr{P}_{\h_*}(O)$. The {\it $H_*$-stability class} $\mathscr{U}_O(H_*)$ is defined by 
$
\mathscr{U}_O(H_*)=\{H\in \mathscr{P}_{\h_*}(O)\, |\, H_*\to H\}
$.
}
\end{define}

\begin{Thm}\label{Coro3}
For every Hamiltonian $H\in \mathscr{P}_{\h_*}(O)$ in the $H_* $-stability class, we have $\mu(H)=\mu(H_*)$.
\end{Thm} 
{\it Proof.} The theorem immediately follows from Proposition \ref{mumu}. $\Box$

\subsection{Basic properties of $\mathscr{U}_O(H_*)$}

In this subsection, we will prove two basic properties of  
$\mathscr{U}_O(H)$.

\begin{Thm}\label{Richness}
For each $H\in \mathscr{P}_{\h_*}(O)$,  the cardinality of $\mathscr{U}_O(H)$ is greater than $\aleph_0$,
the cardinality of the natural numbers. In this sense, $\mathscr{U}_O(H)$ is rich.
\end{Thm}
{\it Proof.}
Suppose that $H$ acts in the Hilbert space $\h$. Note that because $H\in \mathscr{P}_{\h_*}(O)$,
there is a self-dual cone $\Cone$ in $\h$ such that $H\in \mathscr{A}_{\Cone}^+$.
We consider an extended Hilbert space $\h \otimes \BbbC^2$.
We define  a Hamiltonian $H_1$ acting in $\h \otimes \BbbC^2$ by 
$
H_1=H\otimes 1-1\otimes \sigma_1,
$
where $\sigma_1$ is the standard Pauli matrix: $
\sigma_1=\begin{pmatrix}
0 & 1\\
1 & 0
\end{pmatrix}
$. Remark the following facts: 
$\BbbR^2_+=\bigg\{
\left(
    \begin{array}{c}
    x\\
    y
    \end{array}
    \right)
\in \BbbC^2\, \bigg|\, x, y\ge 0\bigg\}$ is a self-dual cone in $\BbbC^2$, and the lowest eigenvalue  of $-\sigma_1$
 is simple with strictly positive eigenvector. Indeed, the eigenvector is given by $\psi_{-\sigma_1}=
 \left(
    \begin{array}{c}
    1/\sqrt{2}\\
    1/\sqrt{2}
    \end{array}
    \right)
 $, which is obviously strictly positive w.r.t. $\BbbR_+^2$.
 Now we define a self-dual cone in $\h \otimes \BbbC^2$ by 
 $
 \Cone_1=\bigg\{\Psi_1\otimes  \left(
    \begin{array}{c}
    1\\
    0
    \end{array}
    \right)
    +\Psi_2\otimes 
    \left(
    \begin{array}{c}
   0\\
    1
    \end{array}
    \right)
    \, \bigg|\, \Psi_1, \Psi_2\in \Cone
    \bigg\}.
 $
 Note that the ground state of $H_1$ is unique and concretely given by $\psi_{H_1}=\psi_H\otimes \psi_{-\sigma_1}$.
Since $H\in \mathscr{A}_{\Cone}^+$,  it holds  that $\psi_H> 0$ w.r.t. $\Cone$. Thus, we readily confirm
 that $\la \Phi|\psi_{H_1}\ra>0$ for all $\Phi\in \Cone_1\backslash \{0\}$, which implies that $\psi_{H_1}>0$
 w.r.t. $\Cone_1$.
 By applying Theorem \ref{Uni}, we conlude that $(H_1+s)^{-1} \rhd 0$ w.r.t. $\Cone_1$ for all $s>-E(H_1)$.
 
We introduce an orthogonal projection $P$ by $P\Psi\otimes r= \Psi\otimes \left(
    \begin{array}{c}
   0\\
    r_2
    \end{array}
    \right)$ for each $\Psi\in \h$ and $r=\left(
    \begin{array}{c}
   r_1\\
    r_2
    \end{array}
    \right) \in \BbbC^2$.
We can identify $\R(P) $ with $\h$
by the isometry $\tau: \R(P)\ni \Psi\otimes 
\left(
    \begin{array}{c}
   0\\
    r_2
    \end{array}
    \right)
\mapsto r_2\Psi \in \h$.
By definition, we have $P\unrhd 0$ w.r.t. $\Cone_1$ and $P \Cone_1=\Cone$ by the aforementioned identification.
Hence, we  readily  check that $H \to H_1$.

Next, let us consider a further extended Hilbert space 
$(\h\otimes \BbbC^2)\otimes \BbbC^2$.
Define a Hamiltonian $H_2$ by 
$
H_2=H_1\otimes 1-1\otimes \sigma_1
$, and define a self-dual cone $\Cone_2$  by $
\Cone_2=
\bigg\{\Phi_1\otimes  \left(
    \begin{array}{c}
    1\\
    0
    \end{array}
    \right)
    +\Phi_2\otimes 
    \left(
    \begin{array}{c}
   0\\
    1
    \end{array}
    \right)
    \, \bigg|\, \Phi_1, \Phi_2\in \Cone_1
    \bigg\}
$.   Using  arguments similar to those in the  last paragraph,  we can  confirm that $
H_1\to H_2
$, which implies that $H\to H_2$.
Repeating this procedure, we can construct a sequence of Hamiltonians $\{H_{\ell}\}_{\ell=1}^{\infty}$ such that $H\to H_{\ell}$.
Therefore, $\mathscr{U}_O(H)$ contains at least countably infinite number of Hamiltonians. $\Box$

\begin{Prop}
Let ${\boldsymbol U}_O$ be the set of all stability classes: ${\boldsymbol U}_O
=\{\mathscr{U}_O(H)\, |\, H\in \mathscr{P}_{\h_*}(O)\}$.
Then ${\boldsymbol U}_O$ is a partially ordered set under set inclusion. 
In addition, the map $\mathscr{U}_O:\, \mathscr{P}_{\h_*}(O) \to {\boldsymbol U}_O$  is 
monotonically decreasing, that is, if $H_1\to H_2$, then $\mathscr{U}_O(H_1) \supseteq \mathscr{U}_O(H_2)$.
\end{Prop}
{\it Proof.} Suppose that $H\in \mathscr{U}_O(H_2)$. Then we have $H_2\to H$. Because $H_1\to H_2$, 
we conclude that $H_1\to H$ by Definition  \ref{HEquiv}. Thus, $H\in \mathscr{U}_O(H_1)$, i.e., $\mathscr{U}_O(H_1) \supseteq \mathscr{U}_O(H_2)$. $\Box$

\subsection{Equivalent  Hamiltonians}

Let $\mathfrak{M}_*$ and $\mathfrak{N}$ be von Neumann algebras on a separable Hilbert spaces $\mathfrak{H}_*$ and $\mathfrak{X}$, respectively.
Let $\Omega_*$ and $\Omega_{\mathfrak{X}}$
 be cyclic and separating vectors  for $\mathfrak{M}_*$ and $\mathfrak{N}$, respectively.
 We denote by $\Cone_*$ and $\Cone_{\mathfrak{X}}$ the natural cones associated with the pairs 
 $\{\mathfrak{M}_*, \Omega_*\}$ and $\{\mathfrak{N}, \Omega_{\mathfrak{X}}\}$, respectively.\footnote{
 We use  $\Delta$ and $J$ to denote  the modular operator and the modular conjugation associated with the pair $\{\mathfrak{M}, \Omega\}$.
  Let 
$
\mathcal{P}_0(\mathfrak{M})=\{A JAJ\, |\, A\in \mathfrak{M}\} .
$
The {\it natural cone} $\Cone$ associated with the pair $\{\mathfrak{M}, \Omega\}$ is defined by 
$
\Cone=\overline{\mathcal{P}_0(\mathfrak{M}) \Omega},
$
where the bar denotes the strong closure. 
It is well-known that  $\Cone$ is a self-dual cone in $\h$ \cite{BR1}.
 }
We set 
\begin{align}
\mathfrak{M}=\mathfrak{M}_*\otimes \mathfrak{N},\ \ \mathfrak{H}=\mathfrak{H}_*\otimes \mathfrak{X},\ \ \Omega=\Omega_*\otimes \Omega_{\mathfrak{X}}. 
\end{align}
Because $\Omega$ is a cyclic and separating vector for $\mathfrak{M}$, we can also define the natural cone $\Cone$ associated with $\mathfrak{M}$.
Physically, $\mathfrak{X}$ represents effects from environments surrounding the system described by $\mathfrak{M}_*$.

Let $\pi=1\otimes |\Omega_{\mathfrak{X}}\ra\la \Omega_{\mathfrak{X}}|$. Then $\pi$ is 
the orthogonal projection from $\mathfrak{H}$ to $\mathfrak{H}_*\cong  \mathfrak{H}_*\otimes \Omega_{\mathfrak{X}}$, where $
\mathfrak{H}_*\otimes \Omega_{\mathfrak{X}}=\{\vphi\otimes \Omega_{\mathfrak{X}}\, |\, \vphi\in \mathfrak{H}_*\}
$. Note that all results in the previous sections  hold true for $\mathfrak{H}_1=\h_*\cong \h_*\otimes \Omega_{\mathfrak{X}}$ and $\mathfrak{H}_2=\h$. 

For given Hamiltonian $H_*$ acting in $\h_*$, 
let us consider $H_*$-stability class $\mathscr{U}_O(H_*)$.

 \begin{define}{\rm 
We say that the Hamiltonian $H$ is  {\it equivalent to $H_*$},  if 
there is a semibounded self-adjoint operator $L$  on $\mathfrak{X}$ 
  such that $L\in \mathscr{A}_{\Cone_{\mathfrak{X}}}^+$ and $H=H_*\otimes 1+1\otimes L$. 
  }
\end{define}

 Let $H\in \mathscr{U}_O(H_*)$. Let $\psi_*$ and $\psi$ be the  normalized ground states of $H_*$ and $H$, respectively.
 The following proposition is readily confirmed.
 \begin{Prop}\label{EquivEasy}
 Suppose that $H$ is equivalent to $H_*$.
 Then there exists a normailzed vector $\omega\in \mathfrak{X}$ such that 
\begin{itemize}
\item[1. ]
$\psi=\psi_*\otimes \omega$;
\item[2. ] $\omega>0$ w.r.t. $\Cone_{\mathfrak{X}}$.
\end{itemize}
 \end{Prop}
 
 \begin{rem}{\rm 
 Let $\{H_{\ell}\}_{\ell=1}^{\infty}$ be the sequence of Hamiltonians constructed in the proof of Theorem \ref{Richness}. By the construction, $H_{\ell+1}$ is equivalent to $H_{\ell}$ for all $\ell\in \BbbN$; thus, every $H_{\ell}$ is equivalent to $H$. In this sense, this example is rivial.
 As we will see in Section \ref{Sec4}, we can construct  sequences of Hamiltonians in $\mathscr{U}_O(H)$ which are inequivalent to each other.
 }
 \end{rem}
 
 As before, we set $\varrho_{\psi_*}=|\psi_*\ra\la \psi_*|$ and $\varrho_{\psi}=|\psi\ra\la \psi|$.
Let $\Tr_{\mathfrak{X}}$ be the partial trace operations  $\Tr_{\mathfrak{X}}: \mathscr{B}(\h) \to \mathscr{B}(\h_*)$. We define a density matrix by $\varrho_{\psi}^{\h_*}=\Tr_{\mathfrak{X}}[\varrho_{\psi}]$.
The quantum relative entropy of $\varrho_{\psi}^{\h_*}$ to $\varrho_{\psi_*}$ is given by 
\begin{align}
S(\varrho_{\psi}^{H_*} | \varrho_{\psi_*})=\Tr[\varrho_{\psi}^{\h_*} \log \varrho_{\psi}^{\h_*}]-\Tr[\varrho_{\psi}^{\h_*} \log \varrho_{\psi_*}].
\end{align}

\begin{define}
{\rm 
We say that $H\in \mathscr{U}_O(H_*)$ is {\it weakly equivalent to  $H_*$}, if $
S(\varrho_{\psi}^{H_*} | \varrho_{\psi_*})=0$.
}
\end{define}

\begin{Thm}\label{TrivialThm}
Let $H\in \mathscr{U}_O(H_*)$. The following (i) and (ii) are equivalent:
\begin{itemize}
\item[(i)] $H$ is weakly equivalent;
\item[(ii)] There exists a normailzed vector $\omega\in \mathfrak{X}$ such that 
\begin{itemize}
\item[1. ]
$\psi=\psi_*\otimes \omega$;
\item[2. ] $\omega>0$ w.r.t. $\Cone_{\mathfrak{X}}$.
\end{itemize}
\end{itemize}
\end{Thm}

\begin{rem}
{\rm

By Proposition \ref{EquivEasy},  the triviality of $H$ implies the condition (ii). However, the converse is  false in general.
}
\end{rem}
{\it Proof of Theorem \ref{TrivialThm}.} (ii) $\Rightarrow$ (i): Easy.

(i) $\Rightarrow$ (ii): First, note that $S(\varrho_{\psi}^{H_*} | \varrho_{\psi_*})=0$ if and only if 
$
\varrho_{\psi}^{H_*}=\varrho_{\psi_*}
$. Thus, there exists an $x\in \mathfrak{X}$ such that $\psi=\psi_*\otimes x$.
Because $\psi$ and $\psi_*$ are strictly positive, $x$ must be strictly positive. $\Box$

\section{Example: construction of a lattice}\label{Sec4}
In this subsection, we will
show that $\mathscr{U}_O(H)$ is truly rich in the  sense that $\mathscr{U}_O(H)$
contains infinitely many {\it inequivalent} Hamiltonians,  and 
  illustrate  the interesting   structure of $\mathscr{U}_O(H)$   by constructing a
specific example.

Let $H_0$ be a self-adjoint operator on $\h_*$, bounded from below.
In this section, we assume the following condition:
\begin{flushleft}
{\bf (H)} $e^{-\beta H_0} \rhd 0$ w.r.t. $\Cone_*$ for all $\beta>0$.
\end{flushleft}
  By Proposition \ref{SemiRe}, we have $H_0\in \mathscr{A}_{\Cone_*}^+$. 
Suppose that $H_0$ commutes with $O$.
Our purpose in this subsection is to examine the stability of $\mu(H_0)$.

For each $n\in \BbbN$ with $n\ge 2$, we consider a Hilbert space 
$\h_*\otimes \BbbC^n$.
Then $\h_*$
 can be regarded as a closed subspace of $\h_*\otimes \BbbC^n$ in the following manner:
 $
 \h_*\cong \h_*\otimes (1/\sqrt{n}, \dots, 1/\sqrt{n})^T\subset  \h_*\otimes \BbbC^n
 $, 
 where $a^T$ indicates the transpose of $a$
  and $\h\otimes a=\{\psi\otimes a| \psi\in \h\}$.
   Thus, $\h_*\otimes \BbbC^n \in \mathbb{H}_{\h_*}$.
 A natural self-dual cone in $\h_*\otimes \BbbC^n$ is  given by 
 \begin{align}
 \Cone_*\otimes \BbbR^n_+=\mathrm{coni}\{ \vphi_j\otimes e_j\, |\, \vphi_j\in \Cone_*,\ j=1, \dots, n\},
 \end{align}
 where $\BbbR_+^n$ is a natural self-dual cone in $\BbbC^n$:  $\BbbR_+^n=\{r=(r_1, \dots, r_n)^T\in \BbbR^n\, |\, r_j\ge 0, j=1, \dots, n\}$, $\{e_j\}_{j=1}^n$ is a standard orthonormal system in $\BbbR^n$ given by $e_j=(0, \dots, \underbrace{1}_{j}, \dots 0)^T$, and $\mathrm{coni}(S)$ is the conical hull of $S$.

Before we proceed, we introduce a useful class of operators.
\begin{define}
{\rm 
Let $\h$ be a Hilbert space and let $\Cone$ be a self-dual cone in $\h$.
 We say that $A\in \mathscr{B}(\h)$ 
is {\it ergodic} w.r.t. $\Cone$ if the following are satisfied:
\begin{itemize}
\item $A\unrhd 0$ w.r.t. $ \Cone$;
\item For each $\xi, \eta\in \Cone \backslash \{0\}$, there exists a $k\in \BbbN\cup \{0\}$ such that $\la \xi|A^k\eta \ra>0$.
\end{itemize}
}
\end{define}

Let $\{n_{\mu}\}_{\mu=1}^{\ell}$ be a set of natural numbers with $\ell\ge 2$ such that $n_1+\cdots+n_{\ell}=N$.
  We set 
\begin{align}
\h_{\mu}= \h_*\otimes \BbbC^{n_{\mu}},\ \ \  \Cone_{\mu}=\Cone_*\otimes \BbbR_+^{n_{\mu}},\ \ \  \mu=1,2, \dots, \ell.
\end{align}
Let $X\in \mathscr{B}(\h_*)$ be self-adjoint. Let $\{Y_{\mu}\}_{\mu=1}^{\ell}$ be a family of self-adjoint operators
such that $Y_{\mu}$ acts in $\BbbC^{n_{\mu}}$. In what follows, we assume the following:
\begin{itemize}
\item[(i)] $X\unrhd 0$ w.r.t. $\Cone_*$;
\item[(ii)] $X$ has purely discrete spectrum and commutes with $O$;
\item[(iii)] $Y_{\mu}$ is ergodic w.r.t. $\BbbR_+^{n_{\mu}}$ for all $\mu=1, \dots, \ell$.
\end{itemize}

\begin{lemm}\label{LemmaY1}
We define a self-adjoint operator $V_{\mu}$ acting in $\h_{\mu}$ by $V_{\mu}=X \otimes Y_{\mu}$.
Then $H_0- V_{\mu}\in \mathscr{A}_{\Cone_{\mu}}^+$. In particular, 
  $H_0- V_{\mu}\in \mathscr{U}_O(H_0)$ for all $\mu=1, \dots, \ell$.
  If $X_{\mu}\neq 1$, then $H_0- V_{\mu}$ is inequivalent to $H_0$. 
\end{lemm}
Note that by the assumptions, $H_0-V_{\mu}$ has purely discrete spectrum.
We will prove Lemma \ref{LemmaY1} in Appendix \ref{PrLemmas}.

Next, let $\mathcal{I}$ be the set of all subsets of $\{1, \dots, {\ell}\}$.
Trivially, $\mathcal{I}$ is a lattice under set inclusion. Let $\mathcal{I}^{\partial}$ be the dual poset of $\mathcal{I}$, that is, the poset with the same underlying set but whose order relation is the opposite of set 
inclusion.
For a given $I=\{{\mu_1}, \dots, {\mu_k}\}$, we set 
\begin{align}
\h_I=\h_* \otimes \BbbC^{n_{\mu_1}}\otimes \cdots \otimes  \BbbC^{n_{\mu_k}},\ \ 
\Cone_I=\Cone_{_*}\otimes \BbbR_+^{n_{\mu_1}}\otimes  \cdots \otimes  \BbbR_+^{n_{\mu_k}}
\end{align}
and $V_I=V_{\mu_1}+\cdots+V_{\mu_{\ell}}$.
Needless to say, $\Cone_I$ is defined by $
\Cone_I=\Big(\cdots\big((\Cone_*\otimes \BbbR_+^{n_{\mu_1}})\otimes \BbbR_+^{n_{\mu_{k+1}}}\big)\otimes \cdots \otimes\BbbR_+^{n_{\mu_{k-1}}} \Big)\otimes \BbbR_+^{n_{\mu_k}}
$.
  If $I=\varnothing$, we simply set $\h_I=\h_*,\ \Cone_I=\Cone_*$ and $V_I=0$.
Note that    
 each $V_{\mu_j}$ acts in $\h_I$ in the following manner:
$
V_{\mu_j}=X\otimes ( \one \otimes \cdots \otimes \overbrace{Y_{\mu_j}}^{j-{\mathrm{th}}}\otimes \cdots \otimes \one). 
$
As before, $\h_*$ can be regarded as a closed subspace of $\h_I$: $\h_*\cong \h_*\otimes \omega_I \subset \h_I$, where 
$\omega_I=\omega_{\mu_1} \otimes \cdots \otimes \omega_{\mu_k}\in \BbbC^{n_{\mu_1}}\otimes \cdots \otimes \BbbC^{n_{\mu_k}}$ with 
$\omega_{\mu_i}=(1/\sqrt{n_{\mu_i}},\dots, 1/\sqrt{n_{\mu_i}})^T \in \BbbC^{n_{\mu_i}}$.

\begin{lemm}\label{LemmaY2}
For each $I\in \mathcal{I}$, we set $H_I=H_0-V_I$. Then 
$H_I\in \mathscr{A}_{\Cone_I}^+$. In particular, 
$H_I\in \mathscr{U}_O(H_0)$. If $X\neq 1$, then $H_I$ is inequivalent to $H_0$. 
\end{lemm}
We will provide a proof of Lemma \ref{LemmaY2} in Appendix \ref{PrLemmas}.

Let $I_1, I_2\in \mathcal{I}$. If $I_1\subseteq I_2$, then $\h_{I_1}$ can be regarded as a subspace of $\h_{I_2}$ in the following manner:
For simplicity, we consider the case where $I_1=\{\mu_1, \dots, \mu_k\}$ and $
I_2=I_1\cup \{\mu_{k+1}, \dots, \mu_{k+\ell}\}$.
Let $\tau$ be a linear operator from $\h_{I_1}$ to $\h_{I_2}$ defined by
\begin{align}
\tau\vphi=\vphi\otimes \omega_{\mu_{k+1}}\otimes \cdots \otimes \omega_{\mu_{k+\ell}}. \label{DefTau}
\end{align}
  It is readily checked that $\tau$ is an isometry. By identifying $\h_{I_1}$ with $\tau \h_{I_1}$,
 $\h_{I_1}$ can be regarded as a subspace of $\h_{I_2}$.
 Note that we can extend this argument to general $I_1$ and $I_2$ with $I_1\subseteq I_2$.

\begin{Thm}
The map $H_{\bullet}: \mathcal{I}^{\partial }\ni I\mapsto H_I\in \mathscr{U}_O(H_0)$ is order-preserving, that is,
if $I_1\subseteq I_2$, then $H_{I_1} \to H_{I_2}$. In particular, $P=\{H_I\}_{I\in \mathcal{I}}$ is a lattice.
The greatest element in $P$ is $H_0$, and the smallest element in $P$ is $H_{\{1, \dots, \ell\}}$.
\end{Thm}
{\it Proof.} For simplicity, we consider the case where $I_1=\{\mu_1, \dots, \mu_k\}$ and $
I_2=I_1\cup \{\mu_{k+1}, \dots, \mu_{k+\ell}\}$.
As explained before, $\h_{I_1}$ can be regarded as a closed subspace of $\h_{I_2}$
 by the isometry $\tau$ defined by (\ref{DefTau}).
 Using this identification, we can identify $\Cone_{I_1}$ with $\tau \Cone_{I_1}
 =\{\psi\otimes \omega_{\mu_{k+1}}\otimes \cdots \otimes \omega_{\mu_{k+\ell}}\, |\, 
 \psi\in \Cone_{I_1}\}$.
 Let $\pi_{I_1, I_2}$ be the orthogonal projection from $\h_{I_2}$ to $\h_{I_1}$:
 \begin{align}
 \pi_{I_1, I_2} \psi=
 P_{\omega_{I_2\backslash I_1} }\psi,
 \end{align}
 where $
 P_{\omega_{I_2\backslash I_1}}=1\otimes |\omega_{I_2\backslash I_1}\ra \la \omega_{I_2\backslash I_1}|
 $
  with $
  \omega_{I_2\backslash I_1}=\omega_{\mu_{k+1}}\otimes \cdots \otimes \omega_{\mu_{k+\ell}}
  $.
  We readily confirm that $\Cone_{I_1}=\pi_{I_1, I_2} \Cone_{I_2}$. Combining this with Lemma \ref{LemmaY2},
  we conclude the assertion in the theorem. $\Box$

\begin{example}{\rm 
For $\ell=3$, we get the following Hasse daigram:
\begin{center}
\includegraphics[width=12cm]{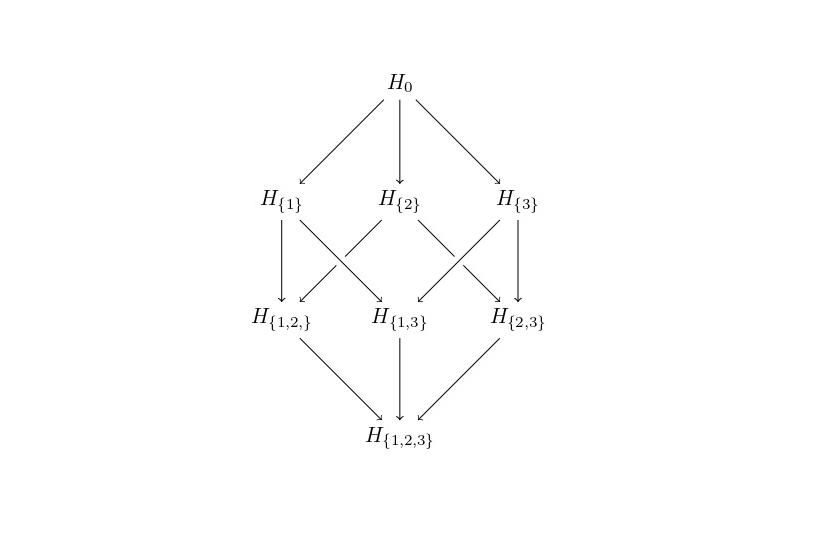}
\end{center}
In the above graph, the vertices are labeld with the elements of the partilly ordered set $P$, and
 the edges indicate the covering relation\footnote{As for the definition of the covering relation, see \cite{Roman}. }.
}
\end{example}

\section{Applications to many-electron systems}
In this section, we will briefly explain  how the theory presented in this paper can be 
applied  to the study of ferromagnetism. Note that the detailed proofs can be 
found in \cite{Miyao8}. 

 Let us consider a finite lattice $\Lambda$. Suppose that $\Lambda$ is bipartite:
 $\Lambda=A\cup B$ with $A\cap B=\varnothing$.
  The Marshall-Lieb-Mattis Hamiltonian is given by 
  \begin{align}
  H_{\mathrm{MLM}}={\boldsymbol S}_A\cdot {\boldsymbol  S}_B,
  \end{align}
  where ${\boldsymbol S}_A=\sum_{x\in A} {\boldsymbol S}_x$ and $
  {\boldsymbol S}_B=\sum_{x\in B} {\boldsymbol S}_x
  $;  ${\boldsymbol S}_x=(S_x^{(1)}, S_x^{(2)}, S_x^{(3)})$ are the spin operators at site
  $x$ satisfying the standard commutation relations:
   \begin{align}
   [S_x^{(j)}, S_y^{(k)}]=i \delta_{xy} \sum_{\ell=1}^3 \varepsilon_{jk\ell} S_x^{(\ell)}.
   \end{align}
   The total spin operators are 
   \begin{align}
   {\boldsymbol S}_{\mathrm{tot}}=\sum_{x\in \Lambda}{\boldsymbol S}_x
   =(S_{\mathrm{tot}}^{(1)}, S_{\mathrm{tot}}^{(2)}, S_{\mathrm{tot}}^{(3)} )
   \end{align}
    and ${\boldsymbol S}_{\mathrm{tot}}^2=\sum_{j=1}^3(S_{\mathrm{tot}}^{(j)})^2$
 with eigenvalues $S(S+1)$.     We say that a  vector $\vphi$ has total spin $S$ if it satisfies $
 {\boldsymbol S}_{\mathrm{tot}}^2\vphi=S(S+1)\vphi
 $.

 We set $\mathfrak{H}_M=\ker(S_{\mathrm{tot}}^{(3)}-M)$, the $M$-subspace.
      We wish to examine properties of $\mathscr{U}_O(H_{\mathrm{MLM}}\restriction \mathfrak{H}_{M=0})$ with 
      $O={\boldsymbol S}_{\mathrm{tot}}^2$.
     The Marshall-Lieb-Mattis theorem \cite{LiebMattis,Marshall} claims that   $\mu(H_{\mathrm{MLM}}
     \restriction \h_0)=
     S_*(S_*+1)$ with $S_*=
     \big||A|-|B|\big|/2$. Hence, we have the following:
        
        \begin{Thm}
        Every Hamiltonian $H$ in $\mathscr{U}_O(H_{\mathrm{MLM}} \restriction \h_0)$ satisfies 
        $\mu(H)=S_*(S_*+1)$.
        \end{Thm}
        \begin{rem}
        {\rm 
        Assume that the ground state of $H$ has total spin $S$.
        We say that the ground state of $H$ exhibits ferromagnetism,
        if $ S$ satisfies $S=c|\Lambda|+o(|\Lambda|)$ with $c>0$ as $|\Lambda|\to \infty$. 
        Therefore, if the ground state of $H_{\mathrm{MLM}} \restriction \h_0$
         exhibits ferromagnetism, then every Hamiltonian in
         $\mathscr{U}_O(H_{\mathrm{MLM}}\restriction \h_0)$ exhibits ferromagnetism as well.
        }
        \end{rem}

       Does $\mathscr{U}_O(H_{\mathrm{MLM}}\restriction \h_0)$ contain physically interesting Hamiltonians? In \cite{Miyao8}, we provide the following answer for this question:
        
        \begin{Thm}\label{StabilityS}
        Let us consider the half-filled many-electron systems on $\Lambda$.
        The following Hamiltonians belong to $\mathscr{U}_O(H_{\mathrm{MLM}} \restriction \h_0)$:
        \begin{itemize}
        \item The Heisenberg Hamiltonian $H_{\mathrm{Heis}}$ restricted to the  $M=0$-subspace;
        \item The Hubbard Hamiltonian $H_{\mathrm{H}}$ restricted to the  $M=0$-subspace;
        \item  The Holstein-Hubbard Hamiltonian $H_{\mathrm{HH}}$ restricted to the  $M=0$-subspace;
        \item A many-electron model coupled to the quantized radiation field $H_{\mathrm{rad}}$ restricted to the  $M=0$-subspace.
        \end{itemize} 
        In addition, these Hamiltonians satisfy the following diagram:
\[
\xymatrix{
& &H_{\mathrm{HH}}\\
H_{\mathrm{MLM}} \equiv H_{\mathrm{Heis}} \ar@{>}[r]& H_{\mathrm{H}}\ar@{>}[ru]\ar@{>}[rd]  &\\
 & &H_{\mathrm{rad}} 
}
\]
(In the above, we abbreviate the restriction of operator $X$  to the $M=0$-subspace, i.e., $X\restriction \ker(S_{\mathrm{tot}}^{(3)})$ to 
$X$. )
        These Hamiltonians  are inequivalent to $H_{\mathrm{MLM}}\restriction \h_0$.
        \end{Thm}
        
        \begin{rem}{\rm 
        \begin{itemize}
        \item Theorem \ref{StabilityS} indicates the stability of Lieb\rq{}s theorem under the influences from environment, e.g., the lattice vibrations and the quantum radiation field.

        \item         
The Marshall-Lieb-Mattis stability class $\mathscr{U}_O(H_{\mathrm{MLM}} \restriction \h_0)$
 is one of the most important examples;
except for this, the Nagaoka-Thouless stability class is  examined in details \cite{Miyao8}.
\end{itemize} 
}
\end{rem}

\appendix
\section{Basic properties of positivity preserving operators}\label{AppA}

\subsection{Positivity preserving operators}

\begin{Prop}\label{Equiv}
Let $A$ be a positive self-adjoint operator. The following statements are equivalent:
\begin{itemize}
\item[{\rm  (i)}] $e^{-\beta A} \unrhd 0$ w.r.t. $\Cone$ for all $\beta \ge 0$.
\item[{\rm (ii)}] $(A+s)^{-1} \unrhd 0$ w.r.t. $\Cone$ for all $s> -E(A)$, where $E(A)=\inf\mathrm{spec}(A)$.
\end{itemize}
\end{Prop}
{\it Proof.} The proposition follows from the following elementary formulas:
\begin{align}
(A+s)^{-1}=\int_0^{\infty}d\beta e^{-\beta (A+s)},\label{ExCh}\\
e^{-\beta A}=\slim \Big(1+\frac{\beta}{n}A\Big)^{-n}.
\end{align}
$\Box$

\begin{Prop}\label{Haagerup}
Let $\Cone$ be a self-dual cone.
Then  $\Cone$ has the following properties:
\begin{itemize}
\item[{\rm (i)}] $\Cone\cap (-\Cone)=\{0\}$.
\item[{\rm (ii)}] There exists a unique antilinear involution $J$ in $\h$ such that
                 $J\xi =\xi$ for all $\xi \in \Cone$.
\item[{\rm (iii)}] Each element $\xi \in \h$ with $J\xi =\xi$ has a unique
                 decomposition $\xi=\xi_+-\xi_-$ where $\xi_+,\xi_-\in\Cone$ and
                 $\la \xi_+| \xi_-\ra=0$.
\item[{\rm (iv)}] $\h$ is linearly spanned by $\Cone$.
\end{itemize}
\end{Prop}
{\it Proof.} See, e.g., \cite[Proof of Proposition 2.5.28 (2), (3) and (4)]{BR1}. $\Box$

\begin{Prop}\label{GSP}
Let $A$ be a positive self-adjoint operator. Assume that $\ex^{-\beta A}
 \unrhd 0$ w.r.t. $\Cone$ for all $\beta \ge 0$. Assume that $E(A)=\inf
 \mathrm{spec}(A)$ is an eigenvalue of $A$. Then there exists a nonzero vector
 $\xi\in \ker(A-E(A))$ such that $\xi\ge 0$ w.r.t. $\Cone$. 
\end{Prop} 
{\it Proof.} {\bf STEP 1.} Let $J$ be an antilinear involution given by
Proposition \ref{Haagerup}.  Set $\mathfrak{H}_J=\{\xi \in \h\, |\, J\xi=\xi\}$.
We will show that $\ker(A-E(A))\cap \h_J\neq \{0\}$.

To see this, let $\xi \in \ker(A-E(A))$. Then we have the decomposition 
$\xi =\Re \xi +i \Im \xi $ with $\Re \xi =\frac{1}{2}(\one +J)\xi $ and $\Im \xi \ 
=\frac{1}{2i}(\one-J)\xi $. Clearly , $\Re \xi, \Im \xi \in \h_J$. Since $\xi \neq 0$, it holds that $\Re \xi \neq 0$
or $\Im  \xi \neq 0$. Since $\ex^{-\beta A}\unrhd 0$ w.r.t. $\Cone$ for all
$\beta \ge 0$, $A$ commutes with $J$.
Thus,  $\Re \xi, \Im \xi\in \ker(A-E(A))\cap \h_J$.

{\bf STEP 2.} Take $\xi \in \ker(A-E(A))\cap \h_J$. By Proposition
\ref{Haagerup} (iii), we have a unique decomposition $\xi=\xi_+-\xi_-$, where 
$\xi_{\pm}\in \Cone$ and $\la \xi_+| \xi_-\ra=0$. Let $|\xi|=\xi_++\xi_-$. Then we
have 
\begin{align}
\ex^{-\beta E(A)}\|\xi\|^2=\la \xi|\ex^{-\beta A}\xi\ra \le \la |\xi||\ex^{-\beta
 A}|\xi|\ra
\le \ex^{-\beta E(A)}\underbrace{\||\xi |\|^2}_{=\|\xi\|^2}.
\end{align} 
Thus,  $|\xi|\in \ker(A-E(A))$. Clearly,  $|\xi|\ge 0$ w.r.t. $\Cone$. $\Box$

\subsection{Positivity improvingness and ergodicity}

\begin{Prop}\label{SemiRe}
Let $A$ be a positive self-adjoint operator.
If $e^{-\beta A} \rhd 0$ w.r.t. $\Cone$ for all $\beta>0$, then 
$(A+s)^{-1} \rhd 0$ w.r.t. $\Cone$ for all $s>-E(A)$.
\end{Prop}
{\it Proof.} This proposition immediately follows from the formula (\ref{ExCh}). $\Box$

\begin{Thm}\label{Uni}
Let $A$ be a positive self-adjoint operator. Suppose that $E(A)=\inf \mathrm{spec}(A)$ is an eigenvalue.
Then the following statements are equivalent:
\begin{itemize}
\item[{\rm (i)}] $(A+s)^{-1} \rhd 0$ w.r.t. $\Cone$ for all $s>-E(A)$.
\item[{\rm (ii)}] $E(A)$ is a simple eigenvalue with a strictly  positive eigenvector w.r.t. $\Cone$. 
\end{itemize}
\end{Thm}
{\it Proof.} This theorem is proved in \cite{Faris}.
Note that the original theorem in \cite{Faris} is  constructed within  real Hilbert spaces, however,
we can readily extend it to a theorem within complex Hilbert spaces. $\Box$

\begin{define}\label{DefOpInq}{\rm 
Let $J$ be the involution given in Proposition \ref{Haagerup}. We set $\h_J=\{\xi\in \h\, |\, J\xi=\xi\}$.
Let $A, B\in \mathscr{B}(\h)$. 
Suppose that $A\h_J\subseteq
 \h_J$ and $B\h_J\subseteq
	     \h_J$. If $(A-B) \Cone\subseteq
	     \Cone$, then we write this as $A \unrhd B$ w.r.t. $\Cone$. 
} 
\end{define}

\begin{Prop}\label{PI1}
Let $A$ be a positive self-adjoint operator and $B$ be a bounded self-adjoint operator.
Suppose that the following conditions are satisfied:
\begin{itemize}
\item[{\rm (i)}] $e^{-\beta A} \unrhd 0$ w.r.t. $\Cone$ for all $\beta \ge 0$;
\item[{\rm (ii)}] $B$ is ergodic w.r.t. $\Cone$.
\end{itemize}
Then $e^{-\beta (A-B)} \rhd 0$ w.r.t. $\Cone$ for all $\beta >0$.
In particular, $(A-B+s)^{-1} \rhd 0$ w.r.t. $\Cone$ for all $s>\inf \mathrm{spec}(A-B)$.
\end{Prop}
{\it Proof.}
By Propositions \ref{ConvexConeOp} and \ref{Equiv}, we find that $
e^{-\beta(A-B)} \unrhd 0
$ w.r.t. $\Cone$ for all $\beta \ge 0$. 
  By the Duhamel formula, we get
\begin{align}
e^{-\beta(A-B)}=\sum_{n=0}^{\infty}\mathcal{I}_n(\beta), \label{DuHa}
\end{align}
where $\mathcal{I}_0(\beta)=e^{-\beta A}$ and 
\begin{align}
\mathcal{I}_n(\beta)=\int_{0\le s_1\le \cdots \le s_n\le \beta} B(s_1)\cdots B(s_n) e^{-\beta A}
\end{align}
with $B(s)=e^{-s A} Be^{sA}$. Note that the right hand side of (\ref{DuHa}) converges in the operator norm topology.
Because $B(s_1)\cdots B(s_n) e^{-\beta A} \unrhd 0$ w.r.t. $\Cone$, provided that $0\le s_1\le \cdots \le s_n\le \beta$, we see that  $\mathcal{I}_n(\beta) \unrhd 0$ w.r.t. $\Cone$ for all $\beta \ge 0$. Thus, we get, by Definition \ref{DefOpInq},
\begin{align}
e^{-\beta (A-B)} \unrhd \mathcal{I}_n(\beta)\ \mbox{w.r.t. $\Cone$} \label{Lower}
\end{align}
for all $n\in \BbbN\cup \{0\}$ and $\beta \ge 0$.

For each $\xi, \eta\in \Cone\backslash \{0\}$, there exists an $\ell\in \BbbN\cup \{0\}$ such that 
$\la \xi|V^{\ell} e^{-\beta A}\eta\ra>0$ due to the ergodicity of $V$.
On the other hand, by (\ref{Lower}), we obtain that 
\begin{align}
\la \xi|e^{-\beta(A-B)} \eta\ra \ge \la \xi|\mathcal{I}_{\ell}(\beta )\eta\ra. \label{Lower2} 
\end{align}
It suffices to show that the right hand side of (\ref{Lower2}) is strictly positive.
To this end, let $F(s_1, \dots, s_{\ell})=\la \xi|B(s_1)\cdots B(s_{\ell})e^{-\beta A}\eta\ra$.
Note that $F(0, \dots, 0)=\la \xi|V^{\ell} e^{-\beta A}\eta\ra>0$ and 
$F(s_1, \dots, s_{\ell}) \ge 0$, provided that $0\le s_1\le \cdots \le s_{\ell}\le \beta$. Because $F$ is continuous in $s_1, \dots, s_{\ell}$, we have
\begin{align}
\la \xi|\mathcal{I}_{\ell}(\beta) \eta\ra=\int_{0\le s_1\le \cdots \le s_{\ell}\le \beta} F(s_1, \dots, s_{\ell})>0.
\end{align}
Thus, we are done. $\Box$

\subsection{Composition of ergodic operators}
\begin{Prop}\label{GeneEr}
Let $\h$ be a Hilbert space and let $\Cone$ be a self-dual cone in $\h$.
Let $A\in \mathscr{B}(\h)$ and let $B\in \mathscr{B}(\BbbC^n)$.
Suppose that $A$ and $B$ satisfy the following conditions:
\begin{itemize}
\item[{\rm (i)}] $A$ is ergodic w.r.t. $\Cone$.
\item[{\rm (ii)}] $B$ is ergodic w.r.t. $\BbbR_+^n$.
 \end{itemize} 
 Then $A\otimes \one +\one \otimes B$ is ergodic w.r.t. $\Cone\otimes \BbbR_+^n$. 
\end{Prop} 
{\it Proof.}
Set $C=A\otimes \one+\one \otimes B$. Take $\vphi, \psi\in (\Cone\otimes \BbbR_+^n) \backslash \{0\}$,
arbitrarily.
We can express $\vphi$ and $\psi$ as 
$\vphi=\sum_{j=1}^n \xi_j\otimes e_j$ and $\psi=\sum_{j=1}^n\eta_j\otimes e_j$, where $\xi_j, \eta_j\in \Cone$, and  $\{e_j\}$ is the  standard orthonormal system in $\BbbR^n$. 
Because $\vphi\neq 0$ and $\psi\neq 0$, there exist $p, q\in \BbbN\cup \{0\}$ such that $\xi_p\neq 0$ and $\eta_q\neq 0$. Thus, we have
\begin{align}
\vphi\ge \xi_p\otimes e_p,\ \ \psi\ge \eta_q\otimes e_q\ \ \mbox{w.r.t. $\Cone\otimes \BbbR^n_+$}. \label{VecInq}
\end{align}

By the assumptions, there exist $M, N\in \BbbN\cup \{0\}$ such that 
$\la \xi_p|A^M \eta_q\ra>0$ and $\la e_p| B^Ne_q\ra>0$.
By the binomial theorem and Definition \ref{DefOpInq}, 
we have
\begin{align}
C^{M+N}=\sum_{j=1}^{M+N} \binom{M+N }{j}
A^{M+N-j}\otimes B^j
\unrhd 
 \binom{M+N}{N}
A^M\otimes B^N\ \ \mbox{w.r.t. $\Cone\otimes \BbbR_+^n$}. \label{BiNo}
\end{align}
Combining (\ref{VecInq}) and (\ref{BiNo}), we get 
\begin{align}
\la \vphi|C^{M+N} \psi\ra \ge 
\binom{M+N }{N}
\la \xi_p|A^M\eta_p\ra\la e_p| B^Ne_q\ra>0.
\end{align}
Thus, we are done. $\Box$

\section{Proof of Lemmas \ref{LemmaY1} and \ref{LemmaY2}}\label{PrLemmas}\label{AppB}

\subsection{A general proposition}
We begin with a general proposition.

\begin{Prop}\label{GenePr}
Let $H_0, O$ and $X$ be self-adjoint operators acting in $\h_*$ satisfying  the all assumptions 
in Section \ref{Sec4}.
Let $Y\in \mathscr{B}(\BbbC^n)$ be a self-adjoint operator satisfying the following condtion:
\begin{flushleft}
{\bf (A)} $Y$ is ergodic w.r.t. $\BbbR_+^n$.
\end{flushleft}
Then $H=H_0\otimes \one- X\otimes Y \in \mathscr{A}_{\Cone_*\otimes \BbbR_+^n}^+$.
\end{Prop}
{\it Proof.} By the Duhamel formula, we have the norm convergent expansion:
\begin{align}
e^{-\beta H}=\sum_{j=0}^{\infty}
 \mathcal{J}_j(\beta), 
 \end{align}
where $\mathcal{J}_0(\beta)=e^{-\beta H_0} \otimes \one$ and 
\begin{align}
\mathcal{J}_j(\beta)=\int_{0\le s_1\le s_2\le \cdots \le s_j\le \beta} X(s_1)\cdots X(s_j) e^{-\beta H_0}\otimes Y^j
\end{align}
with $X(s)=e^{-s H_0} Xe^{sH_0}$. Because $X(s_1)\cdots X(s_j)e^{-\beta H_0} \unrhd 0$ w.r.t. $\Cone_*$,
provided that $0\le s_1\le \cdots \le s_j\le \beta$, we obtain that 
$\mathcal{J}_j(\beta) \unrhd 0$ w.r.t. $\Cone_*\otimes \BbbR_+^n$. Thus, we get
\begin{align}
e^{-\beta H} \unrhd \mathcal{J}_j(\beta) \ \ \mbox{w.r.t. $\Cone_*\otimes \BbbR_+^n$ for all $j$}.
\end{align}

Choose $\vphi, \psi\in (\Cone_*\otimes \BbbR_+^n)\backslash \{0\}$, arbitrarily.
Using an argument similar to that in the proof of Proposition \ref{GeneEr}, we can find $p, q\in \BbbN$
such that $\vphi\ge \xi_p\otimes e_p$ and $\psi\ge \eta_q\otimes e_q$ w.r.t. $\Cone_*\otimes \BbbR_+^n$ with $\xi_p, \eta_q\in \Cone_*\backslash \{0\}$.
Because $Y$ is ergodic w.r.t. $\BbbR_+^n$, there exists an $\ell\in \BbbN\cup \{0\}$ such that 
$\la e_p|Y^{\ell}e_q\ra>0$. For this $\ell$, we claim that 
\begin{align}
\la \xi_p|X(s_1)\cdots X(s_{\ell})e^{-\beta H_0}\eta_q\ra>0, \label{SPV}
\end{align}
provided that $0<s_1<s_2<\cdots<s_{\ell}<\beta$.
To this end, observe that $Xe^{-(\beta-s_{\ell})H_0} \eta_q\ge 0$  and 
$Xe^{-(\beta-s_{\ell})H_0} \eta_q\neq 0
$  by Lemma \ref{VecP2}. Hence, $X(s_{\ell}) e^{-\beta H_0}\eta_q=e^{-s_{\ell} H_0} (
Xe^{-(\beta-s_{\ell})H_0} \eta_q
)>0$ w.r.t. $\Cone_*$ if $0<s_{\ell}<\beta$.
Repeating this argument, we see that 
$
X(s_1)\cdots X(s_{\ell})e^{-\beta H_0}\eta_q>0
$
w.r.t. $\Cone_*$, provided that $0<s_1<s_2<\cdots<s_{\ell}<\beta$.
Therefore, we conclude (\ref{SPV}). To sum, we obtain that 
\begin{align}
\la \vphi|e^{-\beta H}\psi\ra&\ge \la \xi_p\otimes e_p|\mathcal{J}_{\ell}(\beta) \eta_q\otimes e_q\ra\no
&= \int_{0<s_1<s_2<\cdots<s_{\ell}<\beta}  \la \xi_p|X(s_1)\cdots X(s_{\ell})e^{-\beta H_0}\eta_q\ra
\la e_p|Y^{\ell}e_q\ra \no
&>0.
\end{align}
Thus, we are done. $\Box$

\subsection{Proof of Lemma \ref{LemmaY1} }
By  Proposition \ref{GenePr},
we readily confirm that $H_0-V_{\mu} \in \mathscr{A}_{\Cone_{\mu}}^+$.
Recall the identification $\h_*\cong \h_*\otimes \omega_{\mu} \subset \h_{\mu}$, where  $\omega_{\mu}=(1/\sqrt{n_{\mu}},\dots, 1/\sqrt{n_{\mu}})^T\in \BbbC^{n_{\mu}}$.
 Let $\pi$ be the orthogonal projection from $\h_{\mu}$ to $\h_*$ defined by 
 $\pi \psi\otimes a=\la \omega_{\mu}|a\ra\psi$ for each $\psi\in \h_*$ and $a\in \BbbC^{n_{\mu}}$.
 We readily check that $\pi\Cone_{\mu}=\Cone_*$, which implies that $H_0\to H_0-V_{\mu}$. Hence, $H_0-V_{\mu}\in \mathscr{U}_O(H_0)$.
   $\Box$

\subsection{Proof of Lemma \ref{LemmaY2}}

\begin{lemm}\label{SubL}
For each $I=\{\mu_1, \dots, \mu_k\}\in \mathcal{I}$, let 
\begin{align}
Y_I=\sum_{j=1}^k 1\otimes \cdots \otimes \underbrace{Y_{\mu_j}}_{j^{\mathrm{th}}} \otimes \cdots \otimes 1.
\end{align}
Then,  for each $I\in \mathcal{I}$, the following  holds true:
\begin{itemize}
\item[{\bf (E)}]
$Y_I$ is ergodic w.r.t. $\BbbR_+^{n_{\mu_1}}\otimes \cdots \otimes \BbbR_+^{n_{\mu_k}}$.
\end{itemize}
\end{lemm}
{\it Proof.}
We will prove Lemma \ref{SubL} by induction.

Suppose that {\bf (E)} holds true for every $I\in \mathcal{I}$ with $|I|=k$.
Our goal is to prove {\bf (E)} for every $I\in \mathcal{I}$ with $|I|=k+1$.
For a given $I=\{\mu_1, \dots, \mu_{k+1}\}\in \mathcal{I}$, we set $\tilde{I}=\{\mu_1, \dots, \mu_k\}$. Thus,
$I=\tilde{I}\cup \{\mu_{k+1}\}$ holds.
Corresponding to this, $Y_I$ can be expressed as 
$Y_I=Y_{\tilde{I}}\otimes 1+1\otimes Y_{\mu_{k+1}}$.
Because $Y_{\tilde{I}}$ is ergodic w.r.t. $\BbbR_+^{n_{\mu_1}}\otimes \cdots \otimes \BbbR_+^{n_{\mu_k}}$, 
we can apply Proposition \ref{GeneEr} and conclude that $Y_I$ is ergodic 
w.r.t. $\BbbR_+^{n_{\mu_1}}\otimes \cdots \otimes \BbbR_+^{n_{\mu_{k+1}}}$. $\Box$

\begin{flushleft}
{\it Proof of Lemma \ref{LemmaY2}}
\end{flushleft}
Write $I=\{\mu_1, \dots, \mu_k\}$.
Recall the identification $\h_*\cong \h_*\otimes \omega_I\subset \h_I$.
By Proposition \ref{GenePr} and Lemma \ref{SubL}, we 
see that $H_I\in \mathscr{A}_{\Cone_I}^+$.
Let $\pi_I$ be the orthogonal projection from $\h_I$ to $\h_*$ defined by 
$\pi_I \Psi\otimes b=\la \omega_I|b\ra \Psi$ for every $\Psi\in \h_*$ and $b\in 
\BbbC^{n_{\mu_1}}\otimes \cdots\otimes \BbbC^{n_{\mu_k}}$.
We confirm that  $\pi_I \Cone_I=\Cone_*$ holds, which implies that $H_0\to H_I$.
Thus, $H_I\in \mathscr{U}_O(H_0)$.
  $\Box$

\end{document}